\documentclass{amsproc}
\usepackage{amsthm}
\usepackage{amssymb}
\usepackage{amsbsy}
\usepackage{amscd}
\usepackage[mathscr]{eucal}
\usepackage{nicefrac}

\DeclareSymbolFont{rsfs}{U}{rsfs}{m}{n}
\DeclareSymbolFontAlphabet{\mathrsfs}{rsfs}
\usepackage[all]{xy}
\usepackage{tikz}
\usetikzlibrary{arrows,decorations.pathmorphing,backgrounds,positioning,fit,petri}
\usetikzlibrary{decorations.markings}

\usepackage{verbatim}
\usepackage{version}
\usepackage{color}
\definecolor{deepjunglegreen}{rgb}{0.0, 0.29, 0.29}
\newenvironment{NB}{
\color{red}{\bf NB}. \footnotesize
}{}
\excludeversion{NB}
\excludeversion{NB2}
\makeatletter

\usepackage{url}
\usepackage{ifmtarg}

\usepackage{xr-hyper}
\usepackage[%
unicode,
pdftex,
bookmarks=true,
colorlinks=true,
citecolor=deepjunglegreen,
filecolor=deepjunglegreen,
pdfnewwindow=true]{hyperref}
\usepackage{cleveref}
[http://arxiv.org/pdf/1706.05154.pdf]

\usepackage{xspace}

\AtBeginDocument{%
\makeatletter
\let\oldref=\ref
\renewcommand{\ref}[1]{%
  \def\@mystring{earlier-#1}%
  \@ifundefined{r@\@mystring}{%
    \cref{#1}}{%
    \namecref{earlier-#1}\oldref{earlier-#1}%
  }%
}%
\renewcommand{\eqref}[1]{%
  \def\@mystring{earlier-#1}%
  \@ifundefined{r@\@mystring}{%
    \textup{\tagform@{\oldref{#1}}}}{%
    \text{(\oldref{earlier-#1})}%
  }%
}%
\makeatother
}
\crefname{Theorem}{Theorem\xspace}{Theorems}
\crefformat{Theorem}{Theorem~#2#1#3}
\crefname{section}{\S}{\S\S}
\crefformat{section}{\S#2#1#3} % to remove unnecessary space
\crefmultiformat{section}{\S\S#2#1#3}{, #2#1#3}{, #2#1#3}{, #2#1#3}
\crefname{Lemma}{Lemma\xspace}{Lemmas\xspace}
\crefformat{Lemma}{Lemma~#2#1#3}
\crefname{Proposition}{Proposition\xspace}{Propositions\xspace}
\crefformat{Proposition}{Proposition~#2#1#3}
\crefname{Corollary}{Corollary\xspace}{Corollaries\xspace}
\crefformat{Corollary}{Corollary~#2#1#3}
\crefname{Definition}{Definition}{Definitions}
\crefformat{Definition}{Definition~#2#1#3}
\crefname{Remark}{Remark\xspace}{Remarks\xspace}
\crefformat{Remark}{Remark~#2#1#3}
\crefname{Remarks}{Remark\xspace}{Remarks\xspace}
\crefformat{Remarks}{Remark~#2#1#3}
\crefname{Conjecture}{Conjecture\xspace}{Conjectures\xspace}
\crefformat{Conjecture}{Conjecture~#2#1#3}
\crefname{figure}{Figure\xspace}{Figure\xspace}
\crefformat{figure}{Figure~#2#1#3}

\crefformat{appendix}{\S#2#1#3}
\crefformat{enumi}{#2#1#3}

\crefname{equation}{}{}
\crefformat{equation}{(#2#1#3)}

\usepackage[all]{xy}

\usepackage{stmaryrd}
\SetSymbolFont{stmry}{bold}{U}{stmry}{m}{n}
\usepackage{pigpen} % for pigpenfont used for fiber squares
\usepackage{enumitem}
\usepackage{calc}

\usepackage{textgreek}

\usepackage[whole,autotilde]{bxcjkjatype}

\renewcommand{\thesubsection}{\thesection(\@roman\c@subsection)}
\newcounter{number}
\makeatother
\newtheorem{Theorem}[equation]{Theorem}

\newtheorem{Lemma}[equation]{Lemma}

\theoremstyle{definition}

\theoremstyle{remark}
\newtheorem{Remark}[equation]{Remark}

\numberwithin{equation}{section}
\newcommand{\thmref}[1]{\ref{#1}}
\newcommand{\secref}[1]{\ref{#1}}
\newcommand{\defeq}{\overset{\operatorname{\scriptstyle def.}}{=}}
\newcommand{\CC}{{\mathbb C}}
\newcommand{\ZZ}{{\mathbb Z}}

\newcommand{\GL}{\operatorname{GL}}

\newcommand{\algsl}{\operatorname{\mathfrak{sl}}} % because \sl="slant"

\newcommand{\Hom}{\operatorname{Hom}}

\renewcommand{\MR}[1]{}

\newcommand{\dslash}{/\!\!/}
\newcommand{\bN}{\mathbf N}

\newcommand{\tslash}{/\!\!/\!\!/}
\newcommand{\tslabar}{\mathbin{
\setbox0=\hbox{/\!\!/\!\!/}\rule[0.4\ht0]{\wd0}{.3\dp0}\kern-\wd0\box0}}

\makeatletter
\newcommand{\cA}[1][{}]{%
  \@ifmtarg{#1}%
  {\mathcal A}% if #1 is empty
  {\mathcal A(#1)}% if #1 is not empty
}
\newcommand{\cAh}[1][{}]{%
  \@ifmtarg{#1}%
  {\mathcal A_\hbar}% if #1 is empty
  {\mathcal A_\hbar(#1)}% if #1 is not empty
}
\makeatother

\newcommand{\po}{\ar@{}[dr]|{\text{\pigpenfont R}}}
\newcommand{\pb}{\ar@{}[dr]|{\text{\pigpenfont J}}}
\newcommand{\pp}{\ar@{}[dr]|{\text{\pigpenfont P}}}

\newcommand{\cM}{\mathcal M}

\newcommand{\fM}{\mathfrak M}
\newcommand{\fL}{\mathfrak L}
\newcommand{\fT}{\mathfrak T}
\newcommand{\fA}{\mathfrak A}
\newcommand{\bG}{\mathbf G}

\begin{document}
\title[Geometric Satake for affine Lie algebras]
{Geometric Satake correspondence for affine Kac-Moody Lie algebras of type $A$
  \\
\medskip
\small
$A$型のアファイン・カッツ・ムーディー・リー環における幾何学的佐武対応
}
\author[H.~Nakajima]{Hiraku Nakajima (中島　 啓)}
% \address{Research Institute for Mathematical Sciences,
% Kyoto University, Kyoto 606-8502,
% Japan}
% \email{nakajima@kurims.kyoto-u.ac.jp}
\address{Kavli Institute for the Physics and Mathematics of the Universe (WPI),
  The University of Tokyo,
  5-1-5 Kashiwanoha, Kashiwa, Chiba, 277-8583,
  Japan
}
\email{hiraku.nakajima@ipmu.jp}
\thanks{The research of the author is supported by the World Premier
International Research Center Initiative (WPI Initiative), MEXT,
Japan, and by JSPS Kakenhi Grant Numbers
%22244003, % Moriwaki's A
%23224002, % Fukaya's S
%23340005, % my B
%24224001, % Saito's S
%25220701, % Mukai's S
16H06335. % Moriwaki's S
}
%\subjclass[2000]{Primary 22E47; Secondary 14D20,14F43,81T13}

\begin{abstract}
  This is an informal expository article on geometric Satake
  correspondence for affine Kac-Moody Lie algebras of type $A$ given
  in \cite{2018arXiv181004293N}. We emphasize formal analogies between
  this result and the author's earlier results on geometric approaches
  to the representation via quiver varieties.
\end{abstract}

\maketitle

\setcounter{tocdepth}{2}
%\tableofcontents

\section*{}

At Kinosaki Algebraic Geometry Symposium 2018,
the author gave a talk ``Coulomb branches of $3d$ SUSY gauge theories''
explaining
\begin{enumerate}
\item the provisional mathematical definition of Coulomb
branches of $3$-dimensional $\mathcal N=4$ SUSY gauge theories in
\cite{2015arXiv150303676N,2016arXiv160103586B}, and
\item geometric
Satake correspondence for affine Kac-Moody Lie algebras of type $A$ in
\cite{2018arXiv181004293N}.
\end{enumerate}
Since we already have an expository article for the first part
\cite{2017arXiv170605154N} (see also an earlier Japanese version
\cite{2016arXiv161209014N}), we here review only the second part.
We make this article as a continuation of \cite{2017arXiv170605154N},
as we need to presuppose the first part.  When we refer a section
numbered between $1$ and $8$, it means the corresponding section in
\cite{2017arXiv170605154N}.

\subsection*{Acknowledgments}

The author thanks the organizers of Kinosaki Algebraic Geometry
Symposium 2018 for invitation.

\setcounter{section}{8}

\section{Geometric Satake correspondence for complex reductive groups\\\hfill --- disclaimer}

Our main result \ref{thm:satake} can be regarded as an affine Lie
algebra $\algsl(n)_{\mathrm{aff}}$ version of the geometric Satake
correspondence for a complex reductive group $\bG$
\cite{Lus-ast,1995alg.geom.11007G,Beilinson-Drinfeld,MV2}. Our
approach is close to one due to Mirkovi\'c-Vilonen
\cite{MV2}. Fortunately we have many good expository articles on
\cite{MV2}, and there is no reason to try to produce a worse one.
There is another reason to decide not to explain geometric Satake
correspondence: We avoid to use sheaf theoretic language as much as
possible in the spirit of the talk in Algebraic Geometry Symposium.
If we would explain geometric Satake, this decision makes no sense.

We also omit historical accounts. Interested readers should read
\cite[\S3(viii)]{2016arXiv160403625B} and \cite{fnkl_icm}.

\section{Quiver gauge theories of affine type $A$}

%\subsection{Definitions}

Let $n$ be an integer greater than $1$. We have a similar story for
the case $n=1$, but we omit it for brevity.

Let $V = \bigoplus_{i\in\ZZ/n\ZZ} V_i$,
$W = \bigoplus_{i\in\ZZ/n\ZZ} W_i$ be finite dimensional $\ZZ/n\ZZ$
graded complex vector spaces. We define
\begin{equation*}
  G \defeq \prod_i \GL(V_i), \qquad
  \bN \defeq \prod_i \Hom(V_i,V_{i+1}) \oplus \Hom(W_i, V_i).
\end{equation*}
We regard $\bN$ as a representation of $G$ in a natural way. Therefore
we can apply the definition of the Coulomb branch in \secref{sec:def}
for $G$, $\bN$. We denote it by $\cM(\lambda,\mu)$, where $\lambda$,
$\mu$ are weights of the affine Lie algebra $\algsl(n)_{\mathrm{aff}}$
given by
\begin{equation*}
  \lambda = \sum_i \dim W_i\cdot \Lambda_i, \qquad
  \mu = \lambda - \sum_i \dim V_i\cdot \alpha_i,
\end{equation*}
where $\Lambda_i$, $\alpha_i$ are fundamental weights and simple roots
respectively.
These notation are \emph{not} important at this stage, as one can
identify them with dimension vectors
$\mathbf w = (\dim W_i)_{i\in\ZZ/n\ZZ}$,
$\mathbf v = (\dim V_i)_{i\in\ZZ/n\ZZ}$ respectively. Nevertheless it
will be convenient when we relate Coulomb branches to representation
theory of $\algsl(n)_{\mathrm{aff}}$.
The corresponding Higgs branch, which is the symplectic reduction
$\bN\oplus\bN^*\tslash G = \mu^{-1}(0)\dslash G$ (see \secref{sec:CH}),
is nothing but a quiver variety of affine type $A$. We denote it by
$\fM_0(\lambda,\mu)$. In standard notation for a quiver variety, it is
usually denoted by $\fM_0(\mathbf v, \mathbf w)$.
Taking a character $\chi\colon G\to \CC^\times$ given by the product
of determinants of factors $\GL(V_i)$, we can consider a GIT quotient
$\mu^{-1}(0)\dslash_\chi G$. Let us denote it by
$\fM_\chi(\lambda,\mu)$. It is equipped with a projective morphism
$\pi\colon\fM_\chi(\lambda,\mu)\to \fM_0(\lambda,\mu)$. One see that
$\chi$-stability and $\chi$-semistability are equivalent and
$\fM_\chi(\lambda,\mu)$ is smooth.

\section{Quiver varieties and affine Lie algebras}

\subsection{Integrable highest weight representations}

Recall the following, which was proved for more general quiver:
\begin{Theorem}[\cite{Na-quiver,Na-alg}]\label{thm:quiver}
  \textup{(1)} Let $0$ be the point in $\fM_0(\lambda,\mu)$
  corresponding to $0$ in $\bN\oplus\bN^*$. Then
  $\fL_\chi(\lambda,\mu)\defeq \pi^{-1}(0)$ is a lagrangian subvariety in
  $\fM_\chi(\lambda,\mu)$.% It is homotopic to $\fM_\chi(\lambda,\mu)$.

  \textup{(2)} The direct sum
  \begin{equation*}
    \bigoplus_\mu H_{\operatorname{top}}(\fL_\chi(\lambda,\mu))
  \end{equation*}
  of top degree homology groups of $\fL_\chi(\lambda,\mu)$ over $\mu$
  has a structure of an integrable highest weight representation of an
  affine Lie algebra $\algsl(n)_{\mathrm{aff}}$ of type
  $A_{n-1}^{(1)}$.
\end{Theorem}

Recall $\algsl(n)_{\mathrm{aff}}$ is an infinite dimensional Lie
algebra, defined as a central extension of the finite dimensional
complex simple Lie algebra $\algsl(n)$:
$\algsl(n)\otimes\CC[z,z^{-1}] \oplus \CC K \oplus \CC d$. (Here
$d$ is the so-called degree operator.)
Integrable highest weight representations are natural analog of finite
dimensional irreducible representations of $\algsl(n)$. They are
classified by dominant integral weights (i.e.\ highest weights),
$n$-tuple of nonnegative integers in this case. For the representation
in (2), the highest weight is given by $\lambda$, which was
essentially $(\dim W_i)_{i\in\ZZ/n\ZZ}$ as we mentioned above.
Let us denote the corresponding integrable highest weight
representation by $V(\lambda)$. Its weight $\mu$ subspace is denoted
by $V_\mu(\lambda)$.

Let us note that
$\bigoplus_\mu H_{\operatorname{top}}(\fL_\chi(\lambda,\mu))$ has a
distinguished base given by fundamental classes of irreducible
components of $\fL_\chi(\lambda,\mu)$. It is known that the set
$\bigsqcup_\mu \operatorname{Irr}\fL_\chi(\lambda,\mu)$ of irreducible
components has a structure of Kashiwara crystal isomorphic to one for
crystal base of $\mathbf U_q(\algsl(n)_{\mathrm{aff}})$-version for
$V(\lambda)$. See \cite{KS,Saito,Na-Tensor}. We do not give a presice
explain of its meaning (which requires to recall the definition of
Kashiwara crystal): it vaguely means that there is a recursive way to
parametrize irreducible components of $\fL_\chi(\lambda,\mu)$ with
varying $\mu$ starting from the highest weight vector
$\fL_\chi(\lambda,0) = \{0\}$. In fact, this recursive structure was
used to show that the above representation is of highest weight.

Our goal is to relate the Coulomb branch $\cM(\lambda,\mu)$ to
$\algsl(n)_{\mathrm{aff}}$. In fact, we conjecture that we have
a similar relation for Coulomb branches for general quiver gauge
theories. An affine type $A$ quiver is special, as Coulomb branches
have description as Cherkis bow varieties \cite{2016arXiv160602002N},
which can be defined by symplectic reduction of certain finite
dimensional symplectic manifolds (products of vector spaces and
Slodowy slices of type $A$). This description is crucially used in the
proof of the main result \ref{thm:satake} below, although we hide it
from the statement in order to make the length of this article
reasonable.
We lack tools to establish various expected properties of Coulomb
branches of general quiver gauge theories, such as finiteness of
symplectic leaves, description of leaves, etc.

In order to illustrate a formal similarity between two constructions,
one in \ref{thm:quiver} and another in \ref{thm:satake}, let us
briefly review the definition of $\algsl(n)_{\mathrm{aff}}$
representation, formulated in a slightly different manner from the
original.

We use a presentation of $\algsl(n)_{\mathrm{aff}}$ as a
Kac-Moody Lie algebra: it is generated by elements $e_i$, $f_i$, $h_i$
($i\in\ZZ/n\ZZ$), $d$ with certain defining relations which we omit.
Operators $h_i$, $d$ are given by dimensions of $V_i$, $W_i$ on the
component $H_{\mathrm{top}}(\fL_\chi(\lambda,\mu))$, hence not so interesting.
More precisely they are defined so that
$H_{\mathrm{top}}(\fL_\chi(\lambda,\mu))$ is identified with the
weight space $V_\mu(\lambda)$.

Operators $e_i$, $f_i$ are defined via a correspondence in
$\fM_\chi(\lambda,\mu)\times \fM_\chi(\lambda,\mu-\alpha_i)$. Note that
$\mu-\alpha_i$ corresponds to $V\oplus S_i$, where
$S_i$ is a one dimensional vector space with grading $i\in \ZZ/n\ZZ$.

Let $\chi_i\colon G\to \CC^\times$ be a character, which is the
product of determinants of $\GL(V_j)$ for $j\neq i$.
We consider the corresponding GIT quotients
$\fM_{\chi_i}(\lambda,\mu)$, $\fM_{\chi_i}(\lambda,\mu-\alpha_i)$. The
morphism $\pi$ factors as
$\fM_\chi(\lambda,\mu)\to \fM_{\chi_i}(\lambda,\mu)\to
\fM_0(\lambda,\mu)$. The same is true for $\fM_{\chi_i}(\lambda,\mu)$.
Moreover we have a closed embedding
$\fM_{\chi_i}(\lambda,\mu)
\to\fM_{\chi_i}(\lambda,\mu-\alpha_i)$. This is induced from the
corresponding morphism for $\mu^{-1}(0)$ given by extending the data
for $V,W$ to $V\oplus S_i,W$ by zero on the $S_i$ component.  This
extension does not preserve the stability condition for $\chi$, hence
we do not have a morphism
$\fM_\chi(\lambda,\mu)\to\fM_\chi(\lambda,\mu-\alpha_i)$. But it is
well-defined on $\fM_{\chi_i}(\lambda,\mu)$, as we take categorical
quotients with respect to $\GL(V_i)$ at $i$.

We consider the fiber product
\begin{equation*}
  \fM_\chi(\lambda,\mu)\times_{\fM_{\chi_i}(\lambda,\mu-\alpha_i)}
  \fM_\chi(\lambda,\mu-\alpha_i),
\end{equation*}
where $\fM_\chi(\lambda,\mu)\to \fM_{\chi_i}(\lambda,\mu-\alpha_i)$ is
the composition
$\fM_\chi(\lambda,\mu)\to
\fM_{\chi_i}(\lambda,\mu)\to\fM_{\chi_i}(\lambda,\mu-\alpha_i)$. It is
known that the fiber product is lagrangian in
$\fM_\chi(\lambda,\mu)\times \fM_\chi(\lambda,\mu-\alpha_i)$.

There is a distinguished irreducible component $\mathfrak P_i(\lambda,\mu)$ in
the fiber product: it consists of pairs $x$, $x'$ such that $x$ is a
restriction of $x'$ to $V$ (considered as $\subset V\oplus
S_i$). Since $x'$ is a $\GL(V\oplus S_i)$-orbit of linear maps, not a
linear map, its restriction to $V$ needs to be carefully
defined. But we omit this technical detail.

One can check that projections
$p_1, p_2\colon \mathfrak P_i(\lambda,\mu)\to \fM_\chi(\lambda,\mu)$,
$\fM_\chi(\lambda,\mu-\alpha_i)$ are proper. Hence operators
\begin{equation*}
  H_{\mathrm{top}}(\fL_\chi(\lambda,\mu)) \leftrightarrows H_{\mathrm{top}}(\fL_\chi(\lambda,\mu-\alpha_i))
\end{equation*}
are defined by convolution $p_{2*}p_1^*$ and $p_{1*} p_2^*$. These are
definition of elements $e_i$, $f_i$ up to a suitable sign convention.

\subsection{Tensor product}\label{subsec:tensor1}

In \ref{thm:quiver} irreducible integrable representations of
$\algsl(n)_{\mathrm{aff}}$ are constructed. Let us review the
construction of tensor product representations in \cite{Na-Tensor},
which was motivated by earlier works by Lusztig \cite{MR1714628} and
Varagnolo-Vasserot \cite{VV-std}.

Recall that $(\dim W_i)_{i\in\ZZ/n\ZZ}$ gives the highest weight in
\ref{thm:quiver}. We consider a decomposition $W = W^1\oplus W^2$ of
$\ZZ/n\ZZ$-graded vector spaces, which will be related to the tensor
product of representations with highest weights $\lambda^1$,
$\lambda^2$, which correspond to $(\dim W^1_i)$ and $(\dim
W^2_i)$. The following construction depends on the ordering of $W^1$,
$W^2$, though the resulted tensor product does \emph{not}. It is
because the construction can be upgraded to representations of the
quantum loop algebra
$\mathbf U_q(\mathbf L\algsl(n)_{\mathrm{aff}})$ (toroidal
algebra) whose tensor products depend on the order.

Let us take a one parameter subgroup in
$\nu^\bullet\colon\CC^\times\to \GL(W)$ given by
$\nu^\bullet(t) = \operatorname{id}_{W^1}\oplus t\operatorname{id}_{W^2}$.
It acts on $\fM_0(\lambda,\mu)$, $\fM_\chi(\lambda,\mu)$ naturally.
We define attracting subsets in $\fM_0(\lambda,\mu)$:
\begin{equation}\label{eq:91}
\begin{gathered}
%  \begin{split}
  \fT^{\nu^\bullet}_0(\lambda,\mu) \defeq \left\{ x\in\fM_0(\lambda,\mu) \,\middle|\,
    \text{$\lim_{t\to 0} \nu^\bullet(t)x$ exists}\right\}, \\
  \widetilde{\fT}^{\nu^\bullet}_0(\lambda,\mu) \defeq \left\{ x\in\fM_0(\lambda,\mu) \,\middle|\,
    \lim_{t\to 0} \nu^\bullet(t)x = 0\right\}.
%  \fT^{\nu^\bullet}_\chi(\lambda,\mu) & \defeq \pi^{-1}(\fT^{\nu^\bullet}_0(\lambda,\mu)).
%  \end{split}
\end{gathered}
\end{equation}
Their inverse images under $\pi$ are denoted by $\fT^{\nu^\bullet}_\chi(\lambda,\mu)$,
$\widetilde{\fT}^{\nu^\bullet}_\chi(\lambda,\mu)$ respectively.

\begin{Theorem}[\cite{Na-Tensor}]\label{thm:tensor}
  \textup{(1)} The subvariety $\widetilde{\fT}^{\nu^\bullet}_\chi(\lambda,\mu)$ is
  lagrangian in $\fM_\chi(\lambda,\mu)$.

  \textup{(2)} The direct sum
  \begin{equation*}
    \bigoplus_\mu H_{\mathrm{top}}(\widetilde{\fT}^{\nu^\bullet}_\chi(\lambda,\mu))
  \end{equation*}
  of top degree homology groups of
  $\widetilde{\fT}^{\nu^\bullet}_\chi(\lambda,\mu)$ over $\mu$ has a
  structure of an integrable representation of
  $\algsl(n)_{\mathrm{aff}}$, which is isomorphic to the tensor
  product $V(\lambda^1)\otimes V(\lambda^2)$ of two integrable highest
  weights representations with highest weights $\lambda^1$ and
  $\lambda^2$.
\end{Theorem}

The $\algsl(n)_{\mathrm{aff}}$-module structure is given by the
convolution product as in \ref{thm:quiver}. In order to see that the
representation has a correct `size', let us describe irreducible
components of $\widetilde{\fT}_\chi^{\nu^\bullet}(\lambda,\mu)$.
We first observe that the $\nu^\bullet$-fixed point set
$\fM_\chi(\lambda,\mu)^{\nu^\bullet}$ decomposes as
\begin{equation*}
  \fM_\chi(\lambda,\mu)^{\nu^\bullet} \cong \bigsqcup_{\mu=\mu^1+ \mu^2}
  \fM_\chi(\lambda^1,\mu^1)\times \fM_\chi(\lambda^2,\mu^2).
\end{equation*}
The morphism from the right hand side to the left is given by the
direct sum.  To show that it is an isomorphism, one uses that
$\fM_\chi(\lambda,\mu)$ is smooth and a fine moduli
space. See \cite[Lemma~3.2]{Na-Tensor} for detail.
For a point $x\in\widetilde{\fT}^{\nu^\bullet}_\chi(\lambda,\mu)$, the limit
$\lim_{t\to 0} tx$ is a $\nu^\bullet$-fixed point. According to the above
decomposition, we have the induced decomposition
\begin{equation*}
  \widetilde{\fT}^{\nu^\bullet}_\chi(\lambda,\mu) = \bigsqcup_{\mu=\mu^1+ \mu^2}
  \widetilde{\fT}^{\nu^\bullet}_\chi(\lambda^1,\mu^1;\lambda^2,\mu^2).
\end{equation*}
Then $\widetilde{\fT}^{\nu^\bullet}_\chi(\lambda^1,\mu^1;\lambda^2,\mu^2)$ is a
vector bundle over
$\fL_\chi(\lambda^1,\mu^1)\times\fL_\chi(\lambda^2,\mu^2)$, where the
projection to the base is identified with the map given by
$\lim_{t\to 0} tx$.
It consists of points $x$ which are given by `exact sequence'
$0\to x^2 \to x \to x^1\to 0$ with $x^1\in\fL_\chi(\lambda^1,\mu^1)$,
$x^2\in\fL_\chi(\lambda^2,\mu^2)$. See \cite[Remark~3.16]{Na-Tensor}. Now
irreducible components of $\widetilde{\fT}^{\nu^\bullet}_\chi(\lambda,\mu)$ is in
bijection to
\(
\bigsqcup_{\mu=\mu^1+\mu^2} \operatorname{Irr}\fL_\chi(\lambda^1,\mu^1)\times
\operatorname{Irr}\fL_\chi(\lambda^2,\mu^2).
\)

Let us give an example. Let $V = \bigoplus V_i$, $V_i = \CC$
($1\le i\le n-1$), $V_0 = 0$, $W = \bigoplus W_i$, $W_i = \CC$
($i = 1, n-1$), $W_i = 0$ ($i=0$ or $1 < i < n-1$). (We understand
$W_1 = \CC^2$, $W_0 = 0$ when $n=2$.) The corresponding quiver variety
$\fM_0(\lambda,\mu)$ is $\CC^2/(\ZZ/n\ZZ)$, i.e.\ the simple singularity of
type $A_{n-1}$, and $\fM_\chi(\lambda,\mu)$ is its minimal resolution. The
lagrangian subvariety $\fL_\chi(\lambda,\mu)$ consists of a chain of $n-1$ complex
projective lines. We take a decomposition $W = W^1\oplus W^2$ with
$\dim W^1=\dim W^2 = 1$.
The fixed point set $\fM_\chi(\lambda,\mu)^{\nu^\bullet}$
consists of $n$ isolated points, which are south and north poles of
projective lines. There are $(n-2)$ intersection points among them,
and the remaining two are extremal points of the chain.
Thus we have one additional irreducible component in
$\widetilde\fT^{\nu^\bullet}_\chi(\lambda,\mu)$. Depending on a choice of the
decomposition $W=W^1\oplus W^2$, we choose an either of lines through
two extremal points. See~\ref{fig:tensor}. The additional irreducible
component is drawn by a dotted line.

In representation theoretic term this corresponds to the tensor
product representation $\CC^n \times (\CC^n)^*$ of $\algsl(n)$, where
$\CC^n$ is the vector representation and $(\CC^n)^*$ is its dual. It
decomposes into the sum of $\algsl(n)$ (adjoint representation) and
$\CC\operatorname{id}$ (trivial representation). The additional
irreducible component gives the trivial representation summand.

\begin{figure}[htbp]
    \centering
\begin{tikzpicture}[scale=0.3]
\draw (0,10) parabola bend (15,-5) (30,10);
\draw[thick,fill] (9,0) circle (0.3);
\draw[thick,bend left,distance=40] (9,0) to (12,0);
\draw[thick,fill] (12,0) circle (0.3);
\draw[thick,bend left,distance=40] (12,0) to (15,0);
\draw[thick,fill] (15,0) circle (0.3);
\draw[thick,bend left,distance=40] (15,0) to (18,0);
\draw[thick,fill] (18,0) circle (0.3);
\draw[thick,bend left,distance=40] (18,0) to (21,0);
\draw[thick,fill] (21,0) circle (0.3);
\draw[thick,dotted,bend right,distance=25] (21,0) to
node[midway] {
%  $\Leaf_\infty$
}
(28,10);
\end{tikzpicture}
    \caption{$\widetilde\fT^{\nu^\bullet}_\chi(\lambda,\mu)$}
    \label{fig:tensor}
\end{figure}

\begin{Remark}\label{rem:envelope}
  Combining \ref{thm:quiver} with \ref{thm:tensor}, we see that there
  is an isomorphism
  \begin{equation}\label{eq:92}
    H_{\mathrm{top}}(\widetilde{\fT}^{\nu^\bullet}_\chi(\lambda,\mu))
    \cong\bigoplus_{\mu=\mu^1+\mu^2}
    H_{\mathrm{top}}(\fL_\chi(\lambda^1,\mu^1))\otimes
    H_{\mathrm{top}}(\fL_\chi(\lambda^2,\mu^2)).
  \end{equation}
  But the above argument only gives an existence of an isomorphism,
  and does not give a particular isomorphism. Since the tensor product
  $V(\lambda^1)\otimes V(\lambda^2)$ is not irreducible in general,
  such an isomorphism is not unique. For quiver varieties of finite
  type, this ambiguity was fixed by regarding one of tensor factors as
  a lowest weight module. See \cite[Th.~5.9]{Na-Tensor}.
  The stable envelope introduced by Maulik-Okounkov
  \cite{2012arXiv1211.1287M} gives a canonical isomorphism in a
  geometric way for quiver varieties of general types.
\end{Remark}

\begin{Remark}
  The construction can be easily generalized to the case when we
  decompose $W$ into more factors, say $W = W^1\oplus W^2\oplus
  W^3$. It corresponds to the triple tensor product
  $V(\lambda^1)\otimes V(\lambda^2)\otimes V(\lambda^3)$.

  For quiver varieties of affine type $A$ (and more generally when the
  underlying Dynkin diagram contains a loop) there is another type of
  one parameter subgroup acting on $\fM_0(\lambda,\mu)$,
  $\fM_\chi(\lambda,\mu)$: we consider the dilatation action on the
  factor $\Hom(V_{n-1},V_0)$ and the induced action on
  $\fM_0(\lambda,\mu)$, $\fM_\chi(\lambda,\mu)$.
  \begin{NB}
    When the diagram has no loop, this action on $\bN$ can be absorved
    by the $G$-action, hence is trivial on quotient spaces. But this
    is nontrivial when the diagram contains a loop.
  \end{NB}%
  The fixed point set in $\fM_\chi(\lambda,\mu)$ with respect to this
  action is union of quiver varieties of type $A_\infty$, i.e.\ those
  defined by replacing $\ZZ/n\ZZ$-grading above by $\ZZ$-grading. The
  embedding is giving by regarding $\ZZ$-grading modulo $n$.

  The corresponding attracting set realizes a representation which is
  induced by the homomorphism
  $\algsl(n)_{\mathrm{aff}}\to \widehat{\mathfrak{gl}}(\infty)$, where
  the latter is the central extension of the Lie algebra of matrices
  $(a_{ij})_{i,j\in\ZZ}$ with $a_{ij} = 0$ for $|i-j|\gg 0$. (See
  \cite{MR3185361}.)
\end{Remark}

\subsection{Sheaf theoretic formulation}

In order to make analogy between \ref{thm:quiver} and \ref{thm:satake}
closer, we further reformulate the above definition using perverse
sheaves. A reader who is unfamiliar with sheaf theoretical formulation
of convolution products should skip this and next subsections.

It is known \cite[Cor.~10.11]{Na-alg} that
$\pi\colon \fM_\chi(\lambda,\mu)\to \fM_0(\lambda,\mu)$ is semismall
(when we replace the target by the image of $\pi$). Therefore the
direct image $\pi_!(\mathcal C_{\fM_\chi(\lambda,\mu)})$ over
$\fM_\chi(\lambda,\mu)$ is a semisimple perverse sheaf.
Here $\mathcal C_X$ denote the constant sheaf shifted by $\dim X$,
i.e.\ $\CC_X[\dim X]$.
The point $0$ is a stratum and we have
\begin{NB}
\begin{equation*}
  \pi_!(\mathcal C_{\fM_\chi(\lambda,\mu)})\cong
  H_{\mathrm{top}}(\fL_\chi(\lambda,\mu)) \otimes \CC_{0}
  \oplus (\text{other summands}).
\end{equation*}
\end{NB}%
\begin{equation*}
  H_{\mathrm{top}}(\fL_\chi(\lambda,\mu)) \cong
  \Hom(\CC_{0}, \pi_!(\mathcal C_{\fM_\chi(\lambda,\mu)})),
\end{equation*}
where the $\Hom$ in the right hand side is taken in the abelian
category of perverse sheaves on $\fM_0(\lambda,\mu)$ (locally constant along a
certain natural stratification).

As we remarked above, $\pi$ factors through
$\fM_{\chi_i}(\lambda,\mu)$. Let us denote maps by $\pi'$, $\pi''$ so
that $\pi = \pi'\circ\pi''$. Then $\pi''$ is also semismall, and
$\pi''_!(\mathcal C_{\fM_\chi(\lambda,\mu)})$ is a semisimple perverse
sheaf on $\fM_{\chi_i}(\lambda,\mu)$. It decomposes a direct sum of
intersection cohomology complexes of various strata. (One can show
that no intersection cohomology complexes associated with nontrivial
local system appear. cf.\ \cite[Prop.~15.3.2]{Na-qaff}.) Let us write
the decomposition as
\begin{equation*}
  \pi''_!(\mathcal C_{\fM_\chi(\lambda,\mu)})
  \cong \bigoplus_\alpha
  L_\alpha\otimes \mathrm{IC}(\fM_{\chi_i}^\alpha),
\end{equation*}
where $\fM_{\chi_i}^\alpha$ are relevant strata, and $L_\alpha$ is the
top degree homology of the fiber over a point in a stratum
$\fM_{\chi_i}^\alpha$.
It is also known that a stratum is of a form
$\fM_{\chi_i}^{\mathrm{s}}(\lambda,\mu+k\alpha_i)$ with a nonnegative
integer $k$. Here the superscript `s' means the locus of stable
points. See \cite[Prop.~2.25]{Na-branching}.

Since $\pi=\pi'\circ\pi''$, we have
\begin{equation*}
  H_{\mathrm{top}}(\fL_\chi(\lambda,\mu)) \cong
  \bigoplus_\alpha L_\alpha\otimes
  \Hom(\CC_0, \pi'_!\mathrm{IC}(\fM_{\chi_i}^\alpha)).
\end{equation*}
The stratification $\fM_{\chi_i}(\lambda,\mu)=\bigsqcup \fM_{\chi_i}^\alpha$
is compatible with the closed embedding
$\fM_{\chi_i}(\lambda,\mu)\hookrightarrow \fM_{\chi_i}(\lambda,\mu-\alpha_i)$, hence
we may write also
\begin{equation*}
  H_{\mathrm{top}}(\fL_\chi(\lambda,\mu-\alpha_i)) \cong
  \bigoplus_\alpha L_\alpha'\otimes
  \Hom(\CC_0, \pi'_!\mathrm{IC}(\fM_{\chi_i}^\alpha))
\end{equation*}
by setting $L_\alpha = 0$ if $\fM_{\chi_i}^\alpha$ is a new stratum
appearing in $\fM_{\chi_i}(\lambda,\mu-\alpha_i)$.

By the definition of $e_i$, $f_i$ above, they respect the
decomposition, i.e.\ they are tensor products of operators
\(
   L_\alpha \leftrightarrows L_\alpha'
\)
and the identity operator of
\(
   \Hom(\CC_0, \pi'_!\mathrm{IC}(\fM_{\chi_i}^\alpha)).
\)

By a local description of a quiver variety (see \cite[\S3]{Na-qaff}),
the fiber over a point a stratum $\fM_{\chi_i}^\alpha$ is isomorphic
to the central fiber $\fL_\chi$ of another quiver variety, which is of
finite type $A_1$ in the present case.
The construction of $e_i$, $f_i$ is compatible with this description,
namely operators
\(
   L_\alpha \leftrightarrows L_\alpha'
\)
are given by applying the construction in the previous subsection
for quiver varieties of finite type $A_1$. 
Quiver varieties of finite type $A_1$ are cotangent bundles of
Grassmanian \cite[\S7]{Na-quiver}, and the construction coincides with
one done earlier by Ginzburg \cite{MR1111326} in this case.

From this explanation it is clear that dimension of
$\Hom(\CC_0, \pi'_!\mathrm{IC}(\fM_{\chi_i}^\alpha))$ is equal to the
multiplicity of an irreducible representation in
the restriction
\begin{equation*}
  % \bigoplus H_{\mathrm{top}}(\fL_\chi(\lambda,\mu))\!
  V(\lambda)\!
  \downarrow^{\algsl(n)_{\mathrm{aff}}}_{\mathfrak{l}_i}
\end{equation*}
where $\mathfrak l_i$ is the Levi factor corresponding to $i$.
See \cite{Na-branching}.
In the current situation $\mathfrak{l}_i$ is the direct sum of
$\algsl_2$ and an abelian Lie algebra.
The index $\alpha$ corresponds to the highest
weight of an irreducible representation of $\mathfrak{l}_i$. The
nonemptyness of $\fM_{\chi_i}^{\mathrm{s}}(V^0,W)$ implies that
$\dim W_i + \dim V_{i-1}+\dim V_{i+1} \ge 2 \dim V_i$, hence it is
dominant for $\mathfrak l_i$.

\subsection{Sheaf theoretic formulation for tensor products}\label{subsec:sheaf_tensor}

Let us continue sheaf theoretic formulation of results above. We turn
to \ref{thm:tensor}. We review \cite{tensor2}.

We consider the diagram
\begin{equation*}
  \fM_0(\lambda,\mu)^{\nu^\bullet}
  \xleftarrow{p}
  \fT^{\nu^\bullet}_0(\lambda,\mu)
  \xrightarrow{j} \fM_0(\lambda,\mu),
\end{equation*}
where $j$ is the inclusion and $p$ is given by
$\lim_{t\to 0}\nu^\bullet(t)\cdot$. We define the hyperbolic
restriction functor $\Phi$ by $p_* j^!$. Its fundamental properties
are discussed in \cite{Braden,MR3200429}. In particular, $\Phi$ sends
$\pi_!(\mathcal C_{\fM_\chi(\lambda,\mu)})$ to a direct sum of simple
perverse sheaves with shifts.

Moreover, in the current situation,
$\Phi\pi_!(\mathcal C_{\fM_\chi(\lambda,\mu)})$ is a semisimple
perverse sheaf, i.e.\ no shifts are necessary.
This is because $\Phi$ is hyperbolic semi-small in the sense of
\cite{2014arXiv1406.2381B}. This is, in turn, a consequence of
\thmref{thm:tensor}(1).
Now the top degree homology and the hyperbolic restriction are related
by
\begin{equation*}
  H_{\mathrm{top}}(\widetilde{\fT}_\chi^{\nu^\bullet}(\lambda,\mu))
  \cong \Hom(\CC_0,\Phi\pi_!(\mathcal C_{\fM_\chi(\lambda,\mu)})).
\end{equation*}
See \cite[the paragraph after Lemma~4]{tensor2}.

The statement (2) in \ref{thm:tensor} could be regarded as a
computation of the hyperbolic restriction
$\Phi\pi_!(\mathcal C_{\fM_\chi(\lambda,\mu)})$.
Let
\begin{equation*}
  \sigma\colon \bigsqcup_{\mu=\mu^1+\mu^2}
  \fM_0(\lambda^1,\mu^1)\times\fM_0(\lambda^2,\mu^2)
  \to \fM_0(\lambda,\mu)^{\nu^\bullet}
\end{equation*}
denote a morphism given by sum, which is finite and surjective
\cite[Lemma~1]{tensor2}.

\begin{Lemma}[\protect{\cite[Lemma~3]{tensor2}}]
  We have an isomorphism
  \begin{equation*}
    \sigma_! \bigoplus_{\mu=\mu^1+\mu^2}
    (\pi\times\pi)_! \mathcal C_{\fM_\chi(\lambda^1,\mu^1)\times
      \fM_\chi(\lambda^2,\mu^2)} \cong 
    \Phi\pi_!(
    \mathcal C_{\fM_\chi(\lambda,\mu)}).
  \end{equation*}
\end{Lemma}

This is a consequence of description of
$\widetilde{\fT}_\chi^{\nu^\bullet}(\lambda,\mu)$ in \ref{subsec:tensor1}.

Taking a fiber at $0$, we obtain the isomorphism \eqref{eq:92}. As we
mentioned in \ref{rem:envelope}, \eqref{eq:92} is not unique. The same
is true in the above situation. We can also show that the stable
envelope gives a canonical isomorphism also in above. It is a
consequence of \cite[Lemma~4]{tensor2}.

\section{Coulomb branches and affine Lie algebras}

\subsection{Integrable highest weight representations}\label{subsec:integrable2}

We return back to the Coulomb branch $\cM(\lambda,\mu)$ of a quiver
gauge theory of affine type $A$.
As we explained in \ref{sec:def}, the Coulomb branch has a
quantization, and hence a Poisson bracket. It is shown that it gives a
symplectic form on the smooth locus as a statement for general Coulomb
branches. In our case $\cM(\lambda,\mu)$, a finer statement is known:
The decomposition into symplectic leaves are given by
\begin{equation}\label{eq:94}
  \cM(\lambda,\mu) = \bigsqcup_{\kappa,\underline{k}}
  \cM^{\mathrm{s}}(\kappa,\mu)\times S^{\underline{k}}(\CC^2\setminus \{0\}/
  (\ZZ/n\ZZ)),
\end{equation}
where $\underline{k} = [k_1,k_2,\dots]$ is a partition, and $\kappa$
is a dominant weight with
$\lambda-|\underline{k}|\delta\ge\kappa\ge\mu$.  And
$\cM^{\mathrm{s}}(\kappa,\mu)$ is the regular locus of
$\cM(\kappa,\mu)$ if $\dim W\neq 1$ or $\kappa=\mu$ and $\emptyset$
otherwise. See \cite[Th.~7.26]{2016arXiv160602002N}.

Recall that we have chosen a character $\chi\colon G\to \CC^\times$
given by the product of determinants. We consider the induced
homomorphism $\pi_1(\chi)\colon \pi_1(G)\to \pi_1(\CC^\times)$ and its
Pontryagin dual
$\pi_1(\chi)^\wedge \colon \pi_1(\CC^\times)^\wedge \to
\pi_1(G)^\wedge$.
By \ref{subsec:hamtori} $\pi_1(G)^\wedge$ acts on $\cM(\lambda,\mu)$.
In our case we have $\pi_1(G)^\wedge = (\CC^\times)^n$. (We assume
$V_i\neq 0$ for all $i$ for brevity.)
Via $\pi_1(\CC^\times)^\wedge = \CC^\times$, we consider
$\pi_1(\chi)^\wedge$ as a one parameter subgroup in
$\pi_1(G)^\wedge = (\CC^\times)^n$. Let us write $\pi_1(\chi)^\wedge$
by $\chi$ for brevity.

The following was proved in \cite[Prop.~7.30]{2016arXiv160602002N}
(see also \cite[Prop.~4.1]{2018arXiv181004293N}):
\begin{Lemma}\label{lem:fixed}
  The fixed point set $\cM(\lambda,\mu)^{\chi}$ is either empty
  or a single point.
\end{Lemma}

In analogy with \eqref{eq:91} we introduce the attracting set
\begin{equation*}
  \fA_\chi(\lambda,\mu) \defeq
  \left\{ x\in \cM(\lambda,\mu) \,\middle|\,
    \text{$\lim_{t\to 0} \chi(t) x$ exists}\right\}.
\end{equation*}
Note that a tilde version is same by the above lemma.

\begin{Remark}
For Higgs branch, $\chi$ gave a resolution
$\fM_\chi(\lambda,\mu)\to\fM_0(\lambda,\mu)$. The same $\chi$ give a
one parameter subgroup acting on the Coulomb branch
$\cM(\lambda,\mu)$. In the next subsection we see the opposite:
$\nu^\bullet$ acts on $\fM_0(\lambda,\mu)$, and gives a partial
resolution (we will consider the deformation instead) of
$\cM(\lambda,\mu)$.

In physics it is said that the role of FI and mass parameters on
Coulomb and Higgs branches are exchanged. Our definition of Coulomb
branches is given so that this is established rigorously.
\end{Remark}

Now we state the main result
\begin{Theorem}[\cite{2018arXiv181004293N}]\label{thm:satake}
  \textup{(1)} Every intersection with $\fA_\chi(\lambda,\mu)$ and
  a symplectic leaf of $\cM(\lambda,\mu)$ is either empty or a
  lagrangian subvariety in the leaf.

  \textup{(2)} The direct sum
  \begin{equation*}
    \bigoplus_\mu H_{\operatorname{top}}(\fA_\chi(\lambda,\mu))
  \end{equation*}
  of top degree homology groups of $\fA_\chi(\lambda,\mu)$ over $\mu$
  has a structure of an integrable highest weight representation
  $V(\lambda)$ of an affine Lie algebra $\algsl(n)_{\mathrm{aff}}$
  with highest weight $\lambda$.
\end{Theorem}

As we already mentioned, the construction in \ref{thm:satake} is
\emph{formally} similar to one in \ref{thm:quiver}. We define
operators $d$ and $e_i$, $f_i$, $h_i$ ($i\in\ZZ/n\ZZ$). The operators
$d$ and $h_i$ are determined by dimension vectors, namely they are
defined so that $H_{\mathrm{top}}(\fA_\chi(\lambda,\mu))$ is
identified with $V_\mu(\lambda)$.
Operators $e_i$, $f_i$ are defined via study of varieties associated
with another character $\chi_i$. Namely we consider the $\chi_i$-fixed
point set
\(
  \cM(\lambda,\mu)^{\chi_i}
\)
and the attracting set
\(
\fA_{\chi_i}(\lambda,\mu) \defeq \left\{ x\in\cM(\lambda,\mu)
  \,\middle|\, \text{$\lim_{t\to 0}\chi_i(t)x$ exists}\right\}.
\)
We also consider the $\chi_i$-fixed point set
\(
   \fA_{\chi}(\lambda,\mu)^{\chi_i}
\)
in the attracting set $\fA_\chi(\lambda,\mu)$.
Then the original attracting set is the fiber product:
\begin{equation}\label{eq:93}
  \fA_\chi(\lambda,\mu) = \fA_\chi(\lambda,\mu)^{\chi_i}
  \times_{\cM(\lambda,\mu)^{\chi_i}} \fA_{\chi_i}(\lambda,\mu),
\end{equation}
where the morphism
$\fA_{\chi_i}(\lambda,\mu)\to \cM(\lambda,\mu)^{\chi_i}$ is given by
$\lim_{t\to 0}\chi(t)\cdot$.

We then show that $\cM(\lambda,\mu)^\chi$ is either empty or
isomorphic to the Coulomb branch $\cM_{A_1}(\lambda',\mu')$ of type
$A_1$ quiver gauge theory with some $\lambda'$, $\mu'$. Thus
$\fA_\chi(\lambda,\mu)^{\chi_i}$ is the attracting set in the Coulomb
branch of type $A_1$ quiver gauge theory.

This picture enables us to define operators $e_i$, $f_i$ by a
reduction to the $\algsl(2)$ case. This is formally similar to the
method used in \ref{thm:quiver}. A further detail of the definition of
$e_i$, $f_i$ requires hyperbolic restriction functors, and hence is
postponed to \ref{subsec:sheaf}.

Let us also note that
$\bigoplus_\mu H_{\operatorname{top}}(\fA_\chi(\lambda,\mu))$ has a
distinguished base given by fundamental classes of irreducible
components of $\fA_\chi(\lambda,\mu)$. It can be shown that the set
$\bigsqcup \operatorname{Irr} \fA_\chi(\lambda,\mu)$ of irreducible
components has a structure of Kashiwara crystal isomorphic to one for
crystal base of $\mathbf U_q(\algsl(n)_{\mathrm{aff}})$-version of
$V(\lambda)$. See \cite[\S5(vi)]{2018arXiv181004293N}.

\subsection{Tensor product}

Recall that we take a one parameter subgroup
$\nu^\bullet\colon\CC^\times\to\GL(W)$ in \ref{subsec:tensor1}
associated with a decomposition $W = W^1\oplus W^2$. It can be
consider as a flavor symmetry of the gauge theory. As is explained in
\ref{subsec:flavor}, we have the deformation of $\cM(\lambda,\mu)$
parametrized by
$\CC = \operatorname{Spec} H^*_{\CC^\times}(\mathrm{pt})$, as well as
a partial resolution.
Let us use the former and denote its fiber over $1$ by
$\cM^{\nu^{\bullet,\CC}}(\lambda,\mu)$. (The superscript `$\CC$'
indicates that we choose the former. The notation is the same as in
\cite{2018arXiv181004293N}.)

Let $\lambda^1$, $\lambda^2$ denote weights corresponding to
$(\dim W^1)$, $(\dim W^2)$ as before.
As in \ref{lem:fixed} we have
\begin{Lemma}
  The fixed point set $\cM^{\nu^{\bullet,\CC}}(\lambda,\mu)^\chi$ is
  finite, and corresponds naturally to
  \(
  \bigsqcup_{\mu=\mu^1+\mu^2} \cM(\lambda^1,\mu^1)^\chi
  \times \cM(\lambda^2,\mu^2)^\chi.
  \)
\end{Lemma}

Recall $\cM(\lambda,\mu)^\chi$ is either empty or a single point by
\ref{lem:fixed}. And in view of \ref{thm:satake}, it is a single point
if and only if $V_\mu(\lambda)\neq 0$. Therefore the above lemma says
that $\cM^{\nu^{\bullet,\CC}}(\lambda,\mu)^\chi$ is nonempty if and
only if $\mu$ is a weight of $V(\lambda^1)\otimes V(\lambda^2)$. Let
us introduce the attracting set:
\begin{equation*}
  \fA_\chi^{\nu^{\bullet,\CC}}(\lambda,\mu) \defeq
  \left\{ x\in \cM^{\nu^{\bullet,\CC}}(\lambda,\mu) \,\middle|\,
    \text{$\lim_{t\to 0} \chi(t) x$ exists}\right\}.
\end{equation*}

\begin{Theorem}[\protect{cf.~\cite[Cor.~4.9]{2018arXiv181004293N}}]
  The direct sum
  \begin{equation*}
    \bigoplus_\mu H_{\mathrm{top}}(\fA_\chi^{\nu^{\bullet,\CC}}(\lambda,\mu))
  \end{equation*}
  of top degree homology groups of
  $\fA_\chi^{\nu^{\bullet,\CC}}(\lambda,\mu)$ over $\mu$ has a
  structure of an integrable representation of
  $\algsl(n)_{\mathrm{aff}}$, isomorphic to the tensor product
  $V(\lambda^1)\otimes V(\lambda^2)$.
\end{Theorem}

\subsection{Sheaf theoretic formulation}\label{subsec:sheaf}

As in \ref{subsec:sheaf_tensor}, we consider the diagram
\begin{equation*}
  \cM(\lambda,\mu)^{\chi}
  \xleftarrow{p}
  \fA_\chi(\lambda,\mu)
  \xrightarrow{j} \cM(\lambda,\mu),
\end{equation*}
and define the hyperbolic restriction functor $\Phi = p_* j^!$. We
understand $\Phi = 0$ when $\cM(\lambda,\mu)^\chi = \emptyset$. Then
$\Phi(\mathrm{IC}(\cM(\lambda,\mu))$, if it is nonzero, is a
semisimple perverse sheaf over a point $\cM(\lambda,\mu)^\chi$ as
in \ref{subsec:sheaf_tensor}. We have
\begin{equation*}
  H_{\mathrm{top}}(\fA_\chi(\lambda,\mu)) \cong
  \Phi(\mathrm{IC}(\cM(\lambda,\mu))).
\end{equation*}

Let us consider $\chi_i$ as in \ref{subsec:integrable2}. Then we have
hyperbolic restriction functors $\Phi_i$ and $\Phi^i$ associated with
$\fA_\chi(\lambda,\mu)$ and $\fA_\chi(\lambda,\mu)^{\chi_i}$
respectively. The fiber product property \eqref{eq:93} implies
$\Phi = \Phi^i\circ \Phi_i$. Thus
\( \Phi(\mathrm{IC}(\cM(\lambda,\mu))) =
\Phi^i(\Phi_i(\mathrm{IC}(\cM(\lambda,\mu)))).  \)

Recall that $\cM(\lambda,\mu)^{\chi_i}$ is
$\cM_{A_1}(\lambda',\mu')$. The symplectic leaves \eqref{eq:94} for
type $A_1$ are just $\cM_{A_1}(\kappa',\mu')$ where positive integers
$\kappa'$ with $\lambda'\ge\kappa'\ge \mu'$. The morphisms $p$, $j$
are compatible with strata, hence we only get perverse sheaves which
are locally constant along strata. By a little more argument one can
show that no nonconstant local systems appear. Hence
\begin{equation*}
  \Phi_i(\mathrm{IC}(\cM(\lambda,\mu)))
  \cong \bigoplus_{\kappa'} M^{\lambda,\mu}_{\kappa',\mu'}
  \otimes\mathrm{IC}(\cM_{A_1}(\kappa',\mu'))
\end{equation*}
for vector spaces $M^{\lambda,\mu}_{\kappa',\mu'}$.

For $\mu-\alpha_i$, we have
\begin{equation*}
  \Phi_i(\mathrm{IC}(\cM(\lambda,\mu-\alpha_i)))
  \cong \bigoplus_{\kappa'} M^{\lambda,\mu-\alpha_i}_{\kappa',\mu'-2}
  \otimes\mathrm{IC}(\cM_{A_1}(\kappa',\mu'-2)).
\end{equation*}
We then show that there is a natural isomorphism
$M^{\lambda,\mu}_{\kappa',\mu'}\cong
M^{\lambda,\mu-\alpha_i}_{\kappa',\mu'-2}$. This is a consequence of
the factorization property of $\varpi$ introduced in
\ref{subsec:int}. See \cite[Prop.~5.11]{2018arXiv181004293N}.

\bibliographystyle{myamsalpha}
\bibliography{nakajima,mybib,coulomb}    

\def\cprime{$'$} \def\cprime{$'$} \def\cprime{$'$} \def\cprime{$'$}
  \def\cprime{$'$}
  \providecommand{\noopsort}[1]{}\def\cftil#1{\ifmmode\setbox7\hbox{$\accent"5E#1$}\else
  \setbox7\hbox{\accent"5E#1}\penalty 10000\relax\fi\raise 1\ht7
  \hbox{\lower1.15ex\hbox to 1\wd7{\hss\accent"7E\hss}}\penalty 10000
  \hskip-1\wd7\penalty 10000\box7}
\providecommand{\bysame}{\leavevmode\hbox to3em{\hrulefill}\thinspace}
\providecommand{\MR}{\relax\ifhmode\unskip\space\fi MR }
% \MRhref is called by the amsart/book/proc definition of \MR.
\providecommand{\MRhref}[2]{%
  \href{http://www.ams.org/mathscinet-getitem?mr=#1}{#2}
}
\providecommand{\href}[2]{#2}
\begin{thebibliography}{BFN16b}

\bibitem[BD00]{Beilinson-Drinfeld}
A.~Beilinson and V.~Drinfeld, \emph{Quantization of {H}itchin's integrable
  system and {H}ecke eigensheaves}, available at
  \url{http://www.math.uchicago.edu/~mitya/langlands.html}, 2000.

\bibitem[BFN16a]{2016arXiv160103586B}
A.~{Braverman}, M.~{Finkelberg}, and H.~{Nakajima}, \emph{{Towards a
  mathematical definition of Coulomb branches of $3$-dimensional $\mathcal N=4$
  gauge theories, II}}, ArXiv e-prints (2016),
  \href{http://arxiv.org/abs/1601.03586}{{\ttfamily arXiv:1601.03586
  [math.RT]}}.

\bibitem[BFN16b]{2016arXiv160403625B}
\bysame, \emph{{Coulomb branches of $3d$ $\mathcal N=4$ quiver gauge theories
  and slices in the affine Grassmannian (with appendices by Alexander
  Braverman, Michael Finkelberg, Joel Kamnitzer, Ryosuke Kodera, Hiraku
  Nakajima, Ben Webster, and Alex Weekes)}}, ArXiv e-prints (2016),
  \href{http://arxiv.org/abs/1604.03625}{{\ttfamily arXiv:1604.03625
  [math.RT]}}.

\bibitem[BFN16c]{2014arXiv1406.2381B}
\bysame, \emph{{Instanton moduli spaces and $\mathscr W$-algebras}},
  Ast\'erisque (2016), no.~385, vii+128,
  \href{http://arxiv.org/abs/1406.2381}{{\ttfamily arXiv:1406.2381 [math.QA]}}.
  \MR{3592485}

\bibitem[Bra03]{Braden}
T.~Braden, \emph{Hyperbolic localization of intersection cohomology},
  Transform. Groups \textbf{8} (2003), no.~3, 209--216. \MR{1996415
  (2004f:14037)}

\bibitem[DG14]{MR3200429}
V.~Drinfeld and D.~Gaitsgory, \emph{On a theorem of {B}raden}, Transform.
  Groups \textbf{19} (2014), no.~2, 313--358. \MR{3200429}

\bibitem[Fin18]{fnkl_icm}
M.~Finkelberg, \emph{{Double affine {G}rassmannians and {C}oulomb branches of
  $3d$ $\mathcal N=4$ quiver gauge theories}}, Proceedings of the International
  Congress of Mathematicians, 2018 (2017), to appear,
  \href{http://arxiv.org/abs/1712.03039}{{\ttfamily arXiv:1712.03039
  [math.AG]}}.

\bibitem[Gin91]{MR1111326}
V.~Ginzburg, \emph{Lagrangian construction of the enveloping algebra {$U({\rm
  sl}_n)$}}, C. R. Acad. Sci. Paris S\'{e}r. I Math. \textbf{312} (1991),
  no.~12, 907--912. \MR{1111326}

\bibitem[Gin95]{1995alg.geom.11007G}
\bysame, \emph{{Perverse sheaves on a Loop group and Langlands' duality}},
  ArXiv e-prints (1995),
  \href{http://arxiv.org/abs/alg-geom/9511007}{{\ttfamily
  arXiv:alg-geom/9511007 [alg-geom]}}.

\bibitem[KRR13]{MR3185361}
V.~G. Kac, A.~K. Raina, and N.~Rozhkovskaya, \emph{Bombay lectures on highest
  weight representations of infinite dimensional {L}ie algebras}, second ed.,
  Advanced Series in Mathematical Physics, vol.~29, World Scientific Publishing
  Co. Pte. Ltd., Hackensack, NJ, 2013. \MR{3185361}

\bibitem[KS97]{KS}
M.~Kashiwara and Y.~Saito, \emph{Geometric construction of crystal bases}, Duke
  Math. J. \textbf{89} (1997), no.~1, 9--36. \MR{MR1458969 (99e:17025)}

\bibitem[Lus83]{Lus-ast}
G.~Lusztig, \emph{Singularities, character formulas, and a q-analog of weight
  multiplicities}, Ast\'erisque \textbf{101-102} (1983), 208--229.

\bibitem[Lus99]{MR1714628}
\bysame, \emph{Bases in equivariant {$K$}-theory. {II}}, Represent. Theory
  \textbf{3} (1999), 281--353. \MR{1714628}

\bibitem[MO12]{2012arXiv1211.1287M}
D.~{Maulik} and A.~{Okounkov}, \emph{{Quantum Groups and Quantum Cohomology}},
  arXiv e-prints (2012), arXiv:1211.1287,
  \href{http://arxiv.org/abs/1211.1287}{{\ttfamily arXiv:1211.1287 [math.AG]}}.

\bibitem[MV07]{MV2}
I.~Mirkovi{\'c} and K.~Vilonen, \emph{Geometric {L}anglands duality and
  representations of algebraic groups over commutative rings}, Ann. of Math.
  (2) \textbf{166} (2007), no.~1, 95--143. \MR{2342692 (2008m:22027)}

\bibitem[Nak94]{Na-quiver}
H.~Nakajima, \emph{Instantons on {ALE} spaces, quiver varieties, and
  {K}ac-{M}oody algebras}, Duke Math. J. \textbf{76} (1994), no.~2, 365--416.
  \MR{MR1302318 (95i:53051)}

\bibitem[Nak98]{Na-alg}
\bysame, \emph{Quiver varieties and {K}ac-{M}oody algebras}, Duke Math. J.
  \textbf{91} (1998), no.~3, 515--560. \MR{MR1604167 (99b:17033)}

\bibitem[Nak01a]{Na-qaff}
\bysame, \emph{Quiver varieties and finite-dimensional representations of
  quantum affine algebras}, J. Amer. Math. Soc. \textbf{14} (2001), no.~1,
  145--238 (electronic). \MR{MR1808477 (2002i:17023)}

\bibitem[Nak01b]{Na-Tensor}
\bysame, \emph{Quiver varieties and tensor products}, Invent. Math.
  \textbf{146} (2001), no.~2, 399--449. \MR{MR1865400 (2003e:17023)}

\bibitem[Nak09]{Na-branching}
\bysame, \emph{Quiver varieties and branching}, SIGMA Symmetry Integrability
  Geom. Methods Appl. \textbf{5} (2009), Paper 003, 37. \MR{2470410
  (2010f:17034)}

\bibitem[Nak13]{tensor2}
\bysame, \emph{Quiver varieties and tensor products, {I}{I}}, Symmetries,
  Integrable Systems and Representations, Springer Proceedings in Mathematics
  \& Statistics, vol.~40, 2013, pp.~403--428.

\bibitem[Nak16a]{2015arXiv150303676N}
\bysame, \emph{{Towards a mathematical definition of Coulomb branches of
  $3$-dimensional $\mathcal N=4$ gauge theories, I}}, Adv. Theor. Math. Phys.
  \textbf{20} ({\noopsort{2015}}2016), no.~3, 595--669,
  \href{http://arxiv.org/abs/1503.03676}{{\ttfamily arXiv:1503.03676
  [math-ph]}}.

\bibitem[Nak16b]{2016arXiv161209014N}
\bysame, \emph{{Introduction to a provisional mathematical definition of
  Coulomb branches of $3$-dimensional $\mathcal N=4$ gauge theories}},
  第61回代数学シンポジウム報告集 (2016),
  \href{http://arxiv.org/abs/1612.09014}{{\ttfamily arXiv:1612.09014
  [math.RT]}}.

\bibitem[Nak18a]{2017arXiv170605154N}
\bysame, \emph{{Introduction to a provisional mathematical definition of
  Coulomb branches of $3$-dimensional $\mathcal N=4$ gauge theories}}, Modern
  {G}eometry: {A} {C}elebration of the {W}ork of {S}imon {D}onaldson, Proc. of
  Symp. in Pure Math., vol.~99, Amer. Math. Soc., {\noopsort{2017}}2018,
  \href{http://arxiv.org/abs/1706.05154}{{\ttfamily arXiv:1706.05154
  [math.RT]}}, pp.~193--211.

\bibitem[Nak18b]{2018arXiv181004293N}
\bysame, \emph{{Towards geometric Satake correspondence for Kac-Moody algebras
  -- Cherkis bow varieties and affine Lie algebras of type $A$}}, arXiv
  e-prints (2018), arXiv:1810.04293,
  \href{http://arxiv.org/abs/1810.04293}{{\ttfamily arXiv:1810.04293
  [math.RT]}}.

\bibitem[NT17]{2016arXiv160602002N}
H.~Nakajima and Y.~Takayama, \emph{{Cherkis bow varieties and {C}oulomb
  branches of quiver gauge theories of affine type {$A$}}}, Selecta Mathematica
  \textbf{23} (2017), no.~4, 2553--2633,
  \href{http://arxiv.org/abs/1606.02002}{{\ttfamily arXiv:1606.02002
  [math.RT]}}.

\bibitem[Sai02]{Saito}
Y.~Saito, \emph{Crystal bases and quiver varieties}, Math. Ann. \textbf{324}
  (2002), no.~4, 675--688. \MR{MR1942245 (2004a:17023)}

\bibitem[VV02]{VV-std}
M.~Varagnolo and E.~Vasserot, \emph{Standard modules of quantum affine
  algebras}, Duke Math. J. \textbf{111} (2002), no.~3, 509--533. \MR{1885830
  (2003g:17030)}

\end{thebibliography}

\end{document}